\documentclass[12pt,onecolumn,twoside]{IEEEtran}

\usepackage{amsmath,amssymb,euscript,psfrag,latexsym,graphicx}

\begin{document}
\newtheorem{thm}{Theorem}
\newtheorem{cor}[thm]{Corollary}
\newtheorem{lemma}[thm]{Lemma}
\newtheorem{prop}[thm]{Proposition}
\newtheorem{problem}[thm]{Problem}
\newtheorem{remark}[thm]{Remark}
\newtheorem{defn}[thm]{Definition}
\newtheorem{ex}[thm]{Example}
\newcommand{\ignore}[1]{}
\newcommand{\mR}{{\mathbb R}}
\newcommand{\bR}{{\bf R}}
\newcommand{\rR}{{\rm R}}
\newcommand{\mZ}{{\mathbb Z}}
\newcommand{\mC}{{\mathbb C}}
\newcommand{\fR}{{\mathfrak R}}
\newcommand{\fK}{{\mathfrak K}}
\newcommand{\fG}{{\mathfrak G}}
\newcommand{\fT}{{\mathfrak T}}
\newcommand{\fC}{{\mathfrak C}}
\newcommand{\fF}{{\mathfrak F}}
\newcommand{\fM}{{\mathfrak M}}
\newcommand{\mD}{{\mathbb{D}}}
\newcommand{\cU}{{\mathcal{U}}}
\newcommand{\bJ}{{\mathbb{J}}}
\newcommand{\bI}{{\mathbb{I}}}
\newcommand{\s}{{\rm s}}
\newcommand{\ar}{{\rm a}}
\newcommand{\intpi}{{\frac{1}{2\pi}\int_{-\pi}^\pi}}

\title{The maximum entropy ansatz\\
in the absence of a time arrow:\\
fractional-pole models}
\author{Tryphon T. Georgiou\thanks{Department of Electrical and Computer Engineering, University of Minnesota, Minneapolis, MN 55455; {\tt tryphon@ece.umn.edu}} , {\em IEEE Fellow}}
\date{}
\markboth{
\today}
{Georgiou: Randomness and maximum entropy}
\maketitle

\begin{abstract} The maximum entropy ansatz, as it is often invoked in the context of time-series analysis, suggests the selection of a power spectrum which is consistent with autocorrelation data and corresponds to a random process least predictable from past observations. We introduce and compare a class of spectra with the property that the underlying random process is least predictable at any given point from the complete set of past and future observations. In this context, randomness is quantified by the size of the corresponding {\em smoothing error} and deterministic processes are characterized by integrability of the inverse of their power spectral densities---as opposed to the log-integrability in the classical setting. The power spectrum which is consistent with a partial autocorrelation sequence and corresponds to the {\em most random process} in this new sense, is no longer rational but generated by finitely many fractional-poles. 
\end{abstract}


\section{Introduction}
\PARstart{T}{here} is a special place reserved in the spectral analysis literature for the maximum entropy ansatz, and rightly so, due to the multitude of analytic, computational, and practical qualities of maximum entropy spectra. The relevant theory is firmly rooted in analytic interpolation, the moment problem, and the Hilbert space geometry of random processes. The {\em maximum entropy} ({\em ME}) {\em ansatz}, in its basic form, calls for selecting the unique power spectrum which is consistent with a known finite set of autocorrelation moments and is the maximizer of a convex logarithmic functional which represents the entropy rate of the underlying random process. A closely related alternative justification relies on the fact that this {\em maximum entropy process} ({\em ME process}) is the least predictable from past observations and hence, it represents a worst-case situation.

The entropy rate of a random process is an inherently time-dependent concept. This fact becomes apparent in multivariable prediction theory where the variance of optimal least-variance predictors depends on the choice of the time-arrow \cite[Remark 3]{CAR}. It is our contention that often there is no natural direction of time. This is the case when statistics are obtained from an array of sensors and the index of the autocorrelation moments represents spacial separation. It is the case when we consider sparse records with both, past and future data available but with possible gaps. It is also the case, when we want to estimate the power spectrum and have no plans to use it for prediction in one way or another. In all such cases the rationale of the ME ansatz may be called into question. Hence, the purpose of this work is to study a time-arrow independent counterpart. In this, a power spectrum is selected so that the underlying random process evaluated at any point in time is the least predictable from the complete set of all other past and future values. In other words, it is the variance of the optimal smoothing filter which is sought to be maximal, as opposed to the variance of the optimal (time-arrow dependent) predictor.

Power spectra which are consistent with a finite set of (contiguous) autocorrelation statistics and correspond to a worst-case smoothing error for the relevant random process, turn out to have an all-pole representation, very much like the ones that result in from the ME ansatz but with one important difference. These spectra are inverses of the {\em square root} of positive trigonometric polynomials, and hence, their poles are fractional. They also share a similar property with ME spectra in that they are extrema of a corresponding convex functional ---which, however, is no longer logarithmic. Computation of their respective parameters is slightly more involved than having to solve linear (Yule-Walker-Levinson) equations. They can be computed as fixed points of suitable differential equations
originating from a homotopy-based method in determining functional extrema. For convenience, and lacking a better terminology, we refer to this new class of spectra and the respective processes as {\em most random} ({\em MR}).

The maximum entropy ansatz has a fifty year history or more. We will not attempt to overview significant milestones but refer to \cite{StoicaMoses} for textbook exposition of relevant material, to \cite{Haykin,LevineTribus} for an overview of relevant research in signal processing, to Burg \cite{burg} who is credited with introducing the maximum entropy ansatz in time series analysis, and to Jaynes \cite{jaynes} and Csiszar \cite{csiszar2} for systematic analyses of the ansatz and its relevance in scientific modeling.

\section{Development and main results}
\newcommand{\cE}{{\mathcal E}}
\newcommand{\me}{{\rm _{ME}}}
As explained in the introduction, we consider the problem of spectral analysis based on partial autocorrelation statistics. 
Thus, we begin with a finite set of autocorrelation samples $R_k:=\cE\{u_\ell u_{\ell-k}^*\}$, for $k=0,1,2,\ldots,n$, of a zero-mean, stationary scalar random process $\{u_\ell\;:\; \ell\in\mZ\}$, where ``$^*$'' denotes complex conjugation (together with transposition when applied to vectorial quantities). The discrete ``time index'' may represent a spatial coordinate when the $u_\ell$'s are readings at, say, a number of uniformly and linearly spaced sensor locations.

Without loss of generality we assume that
\[
\bR_n:=\left[\begin{array}{cccc}R_0&R^*_1&\ldots&R^*_n\\
                                                 R_1&R_0&\ldots &R^*_{n-1}\\
                                                 \vdots& & \ddots& \vdots\\
                                                 R_n& R_{n-1} &\ldots& R_0\end{array}\right]> 0,
\]
i.e., that it is {\em positive definite},
for otherwise there is a unique power spectrum $d\mu(\theta)$ for which
\begin{equation}\label{moment}
R_k=\frac{1}{2\pi}\int_{-\pi}^{\pi}e^{-jk\theta}d\mu(\theta),\mbox{ for }k=0,1,\ldots,n
\end{equation}
see e.g., \cite{StoicaMoses,GrenanderSzego}.
The following theorem summarizes known facts about the maximum entropy power spectrum which is consistent with $\bR_n$.

\begin{thm}\label{MEtheorem}{\sf
Provided $\bR_n>0$ there exists a unique power spectrum $d\mu_\me$ (i.e., a nonnegative measure on $[-\pi,\,\pi)$) which satisfies (\ref{moment}) and is a maximizer of the following convex functional
\begin{equation}
\bI(d\mu/d\theta):=\frac{1}{2\pi}\int_{-\pi}^\pi \log\left(\left(\frac{d\mu(\theta)}{d\theta}\right)^{-1}\right)d\theta.
\end{equation}
Further, $d\mu_\me$  is absolutely continuous (with respect to the Lebesgue measure) and of the form
\[
d\mu_\me(\theta)=f_\me(\theta)d\theta
\]
where the spectral density $f_\me(\theta)$ is the inverse of a positive trigonometric polynomial of degree at most $n$, i.e.,
\[
f_\me(\theta)=\frac{k_\me^2}{|a(e^{j\theta})|^2}
\]
with $k_\me^2>0$, and $a(z)=1+a_1z+\ldots+a_nz^n$. The polynomial $a(z)$ can be selected to have all of its roots in the complement $\mD^c:=\{z\,:\,|z|> 1\}$ of the unit disc $\mD$ ($\mD:=\{z\,:\,|z|\leq 1\}$) of the complex plane, in which case
\begin{equation}\label{alphas}
\alpha_k=\left\{\begin{array}{c} -a_k\\0\end{array} \begin{array}{l} \mbox{ for }k\leq n\\\mbox{ otherwise}\end{array}\right.
\end{equation}
is the (unique) minimizer of the variance
\[
\cE_{d\mu_\me}\{|u_0 -\hat{u}_{0|{\rm past}}|^2\}
\]
of the (one-step-ahead) prediction error when the {\em predictor}
\begin{equation}\label{upast}
\hat{u}_{0|{\rm past}}:=\sum_{k>0} \alpha_k u_{-k}
\end{equation}
is sought to depend linearly on past observations.
In general, the minimal variance of the prediction error depends on the choice of $d\mu$ (which is subject to (\ref{moment})).
This variance is maximal when $d\mu_\me$ is selected, i.e.,
the maximum entropy power spectrum solves the min-max problem:
\[
\max_{d\mu} \min_{\alpha_k,\;k>0}\left\{ \cE_{d\mu}\{|u_0 -\sum_{k>0} \alpha_k u_{-k}|^2\}\;:\; \mbox{ } \mbox{(\ref{moment}) holds} \right\}.\Box
\]
}\end{thm}

In the theorem and throughout, $d\mu/d\theta =f$ denotes the power spectral density function which is independent of any possible singular part of the spectral measure $d\mu$.
The theorem is well known and has its roots in the classical theory of moments and the theory of orthogonal polynomials. For a proof see \cite{Geronimus,GrenanderSzego}, cf.\ \cite{StoicaMoses}.  More specifically, the extremal properties of $a(z)$ are established in e.g., \cite[page 38]{GrenanderSzego}, see also \cite[Chapter VIII]{Geronimus}. The fact that $f_\me$ is consistent with the autocorrelation moment constraints inherited by $\bR_n$ follows from \cite[Equations (1.15), (1.18)]{Geronimus}. On the other hand, the entropy functional $\bI(\cdot)$ is clearly convex and a variational argument easily shows that the minimizer is of the form indicated. The last statement follows from the fact that (see \cite[page 38, section 2.2]{GrenanderSzego})
\[
\min_{\alpha_k,\;k>0}
\cE_{d\mu_\me}\{|u_0 -\sum_{k>0} \alpha_k u_{-k}|^2\}
=\frac{   \det \bR_n   }{    \det \bR_{n-1}   }
\]
is achieved for the choice $\alpha_k=a_k$, while
\[
\cE_{d\mu}\{|u_0 -\sum_{k>0} \alpha_k u_{-k}|^2\}=\frac{\det\bR_{n}}{\det\bR_{n-1}}
\]
is independent of $d\mu$ as long as (\ref{moment}) holds. An alternative derivation of all the claims in the theorem can be contructed in a way analogous to the steps used in the proof of Theorem \ref{thm2} below, which we provide in Section \ref{proofs}.

The functional $\bI(\cdot)$ in Theorem \ref{MEtheorem} can be interpreted to represent {\em entropy rate} (see \cite{Haykin}) and has been introduced into time-series modeling by Burg \cite{burg}.
It is also interesting to note that the maximum entropy solution $d\mu_\me$ together with $\alpha_k$'s in (\ref{alphas}) represent a saddle point of $\cE_{d\mu}\{|u_0 -\sum_{k>0} \alpha_k u_{-k}|^2\}$ seen as function of two variables, $d\mu$ and the infinite coefficient vector $(\alpha_1,\alpha_2,\ldots\,)$.

An alternative choice for a solution to (\ref{moment}) corresponding to the least predictable process (MR-process) from combined past and future values can be also obtained via convex optimization of a suitable functional. The following proposition presents this MR-solution and highlights its justification as the worst-case senario with regard to a corresponding smoothing problem. The development mirrors the case of the ME-solution. 

\newcommand{\mr}{{\rm _{MR}}}
\begin{thm}\label{thm2}{\sf
Provided $\bR_n>0$ there exists a unique power spectrum $d\mu_\mr$ (nonnegative measure on $[-\pi,\,\pi]$) which satisfies (\ref{moment}) and is a minimizer of the following concave functional
\begin{equation}
\bJ(d\mu/d\theta):=\frac{1}{2\pi}\int_{-\pi}^\pi \left(\frac{d\mu(\theta)}{d\theta}\right)^{-1}d\theta.
\end{equation}
Further, $d\mu_\mr$  is absolutely continuous (with respect to the Lebesgue measure) and of the form
\[
d\mu_\mr(\theta)=f_\mr(\theta)d\theta
\]
where the spectral density $f_\mr(\theta)$ is the square root of the inverse of a positive trigonometric polynomial of degree at most $n$, i.e.,
\[
f_\mr(\theta)=\frac{k_\mr^2}{\sqrt{b(e^{j\theta})}}
\]
with $k_\mr^2>0$ and
\[
b(e^{j\theta})=b_{-n}e^{-nj\theta}+\ldots+b_0+\ldots+b_ne^{nj\theta}>0
\]
for $\theta\in[-\pi,\,\pi]$ (and $b_{-k}:=b^*_k$).  The constant $k_\mr^2$ can be selected so that the  trigonometric polynomial $b(e^{j\theta})$ satisfies
\[\frac{1}{2\pi}\int_{-\pi}^{\pi}\sqrt{b(e^{j\theta})}d\theta = 1,
\]
in which case,
\begin{equation}\label{betas}
\beta_k=\left\{\begin{array}{c} -\rho_k\\ 0\end{array}\begin{array}{l}\mbox{ when }1< |k|\\ 0 \mbox{ when }k=0\end{array}\right.
\end{equation}
with $\rho_\ell$ the coefficients of the Fourier series of
\[\sqrt{b(e^{j\theta})}= \ldots+\rho_{-2}e^{-2j\theta}+\rho_{-1}e^{-j\theta}+1+\rho_{1}e^{j\theta}+\rho_{2}e^{2j\theta}+\ldots
\]
is the (unique) minimizer
of the variance
\[
\cE_{d\mu_\mr}\{|u_0 -\hat{u}_{0|{\rm past \;\&\; future}}|^2\}
\]
of the smoothing error when the {\em smoothing filter}
\begin{equation}\label{uboth}
\hat{u}_{0|{\rm past \;\&\; future}}:=\sum_{k\neq 0} \beta_k u_{-k}
\end{equation}
is sought to depend linearly on past and future observations.
In general, the minimal variance of the error depends on the choice of $d\mu$ (which is subject to (\ref{moment})).
This variance is maximal when $d\mu_\mr$ is selected, i.e.,
the most random power spectrum solves the min-max problem:
\[
\max_{d\mu} \min_{\beta_k,\;k\neq 0}\left\{ \cE_{d\mu}\{|u_0 -\sum_{k\neq 0} \beta_k u_{-k}|^2\}\;:\; \mbox{ } \mbox{(\ref{moment}) holds}\right\}.
\]
}\end{thm}

The last statement of the theorem echoes the analogous property of the maximum entropy solution. In fact, it can be seen that in the present case $d\mu_\mr$, together with the coefficients $\beta_k$'s in (\ref{betas}) for the smoothing filter, represent a saddle point of $\cE_{d\mu}\{|u_0 -\sum_{k\neq 0} \beta_k u_{-k}|^2\}$.

The ME-power spectrum is rational and its coefficients can be obtained by solving a system of linear equations (the Yule-Walker-Levinson equations) which give rise to the following expression for
\begin{equation}
a(z)=\frac{1}{\det(\bR_{n-1})}\det\left({\footnotesize
\begin{array}{cccc}
R_0&R_{-1}&\ldots&R_{-n}\\
R_1&R_0&\ldots &R_{-n+1}\\
\vdots &&&\vdots\\
R_{n-1}&R_{n-1}&\ldots&R_{-1}\\
z^n&z^{n-1}&\ldots& 1
\end{array}
}\right),
\label{orthogonalpoly}
\end{equation}
while $k_\me^2=\det(\bR_n)/\det(\bR_{n-1})$
e.g., see \cite{StoicaMoses} and also \cite[page 156]{Geronimus}.
The corresponding random process can then be simulated via a Markovian realization---in fact via an autoregressive model with transfer function $k_\me/a(z)$ driven by a unit-variance, white-noise input, cf.\ \cite{StoicaMoses}.

The case of the MR-power spectrum differs substantially in this respect. The power spectral density function is not rational. However, its coefficients can be readily obtained from the data $\bR_n$ using the formalism in \cite{ac_may2004,IT}.
This is explained in the following statement.

\begin{thm}\label{computation}{\sf Let $\bR_n>0$, define the column vectors
\begin{eqnarray*}
\rR_1&:=&[\begin{array}{ccccccc}R_n^*&\ldots&R_1^*&R_0&R_1&\ldots&R_n\end{array}]^\prime, \mbox{ and}\\
G(e^{j\theta})&:=&[\begin{array}{ccccccc}e^{jn\theta}&\ldots&e^{j\theta}&1&e^{-j\theta}&\ldots&e^{-jn\theta}\end{array}]^\prime,
\end{eqnarray*}
of size $2n+1$, where $^\prime$ denotes transposition (without complex conjugation),
and let the $\lambda(t)\in\mC^{(2n+1)\times 1}$ represent the solution of the differential equation
\begin{equation}\label{fixedpoint}
\frac{d \lambda(t)}{dt}=M(\lambda(t))^{-1}\left(\rR_1-\intpi \frac{G(e^{j\theta})}{\sqrt{\lambda(t)G(e^{j\theta})}}d\theta\right)
\end{equation}
on $[0,\infty)$, where
\begin{equation}\label{Jacobian}
M(\lambda(t)):=-\intpi G(e^{j\theta})\frac{1/2}{\left(\lambda(t)G(e^{j\theta})\right)^{3/2}}G(e^{j\theta})^*d\theta
\end{equation}
and
\[
\lambda(0)=\lambda_0:=[\;
\underbrace{\begin{array}{ccc}0&\ldots&0\end{array}}_n
1 \underbrace{\begin{array}{ccc}0&\ldots&0\end{array}}_n\;].
\]
Then the following hold:
\begin{itemize}
\item[(i)] $\lambda(t)$ tends to a limit $\lambda_\mr\in\mC^{(n+1)\times 1}$ as $t\to\infty$,
\item[(ii)] $\lambda_\mr G(e^{j\theta})>0$ for all $\theta\in[-\pi,\,\pi)$ and
\begin{equation}\label{conditionii}
d\mu(\theta)=\frac{1}{\sqrt{\lambda_\mr G(e^{j\theta})}}d\theta \mbox{ satisfies }(\ref{moment}),
\end{equation}
\item[(iii)] $\lambda_\mr$ is the unique value in $\mC^{2n+1}$ for which (ii) holds.
\end{itemize}
}
\end{thm}


\section{Notation and preliminaries}

We consider the scalar zero-mean stationary random process
$\{u_k,\;k\in\mZ\}$ and, as before, we let
$R_0,\,R_1,\,R_2,\ldots$ represent its sequence of autocorrelation samples and 
$d\mu(\theta)$ its power spectrum.
We study quadratic optimization problems with respect to the usual inner product
\begin{eqnarray}\nonumber
\langle \sum_k a_ku_k,\sum_\ell b_\ell u_\ell\rangle_{d\mu}&:=&\cE_{d\mu}\{ (\sum_k a_ku_k)(\sum_\ell b_\ell u_\ell)^*\}\\
&=&\sum_{k,\ell} a_kR_{k-\ell}b_\ell^* \label{Tproduct}
\end{eqnarray}
where $R_{-m}:=R_m^*$.
As usual \cite{masani}, the closure of ${\rm span}\{u_k\,:\,k\in\mZ\}$, which we denote by $\cU$, can be identified with the space $L_{2, d\mu}$ of functions which are square integrable with respect to $d\mu(\theta)$ on the unit circle with inner product
\[
\langle a,b\rangle_{d\mu}:=
\frac{1}{2\pi}\int_{-\pi}^\pi a(\theta)(b(\theta))^*d\mu(\theta)
\]
where $a(\theta)=\sum_k a_ke^{jk\theta}$ and $b(\theta)=\sum_\ell b_\ell e^{j\ell\theta}$.
Then it is quite standard that the correspondence
\begin{eqnarray*}
\cU\to L_{2,d\mu}&:& \sum_k a_ku_k \mapsto \sum_k a_ke^{jk\theta}
\end{eqnarray*}
is a Hilbert space isomorphism.

Least-variance approximation problems can equivalently be expressed in $L_{2,d\mu}$.
In particular, the variance $\cE_{d\mu}\{|u_0 -\hat{u}_{0|{\rm past}}|^2\}$ of the one-step-ahead prediction error
\[
u_0-\hat{u}_{0|{\rm past}}
\]
with $\hat{u}_{0|{\rm past}}=\sum_{k>0}\alpha_ku_{-k}$ as in (\ref{upast}), can equivalently be expressed in the form
\begin{equation}\label{infalphas}
\|1-\sum_{k>0}\alpha_ke^{jk\theta}\|_{d\mu}^2,
\end{equation}
and similarly, the variance of the smoothing error $\cE_{d\mu}\{|u_0 -\hat{u}_{0|{\rm past \;\&\; future}}|^2\}$ is simply
\begin{equation}\label{infbetas}
\|1-\sum_{k\neq 0}\beta_ke^{jk\theta}\|_{d\mu}^2
\end{equation}
in view of $\hat{u}_{0|{\rm past \;\&\; future}}$ given in (\ref{uboth}).

The power spectrum $d\mu$ is a bounded nonnegative measure on $[-\pi,\,\pi)$ and admits a decomposition $d\mu=d\mu_\s+fd\theta$ with $d\mu_\s$ a singular measure and $fd\theta$ the absolutely continuous part of $d\mu$ (with respect to the Lebesgue measure). Then, the variance of the optimal
one-step-ahead prediction error is given in terms of the power spectral density function $f$ by the celebrated Szeg\"{o}-Kolmogorov formula given below.

\begin{thm} \cite{Szego} {\sf With $d\mu=d\mu_\s+fd\theta$ as above
\[
\inf_\alpha \|1-\sum_{k>0}\alpha_ke^{jk\theta}\|_{d\mu}^2=\exp\left\{\frac{1}{2\pi}\int_{-\pi}^\pi \log f(\theta)d\theta\right\}
\]
when $\log f\in L_1$, and zero otherwise.
}\end{thm}

For a proof see \cite[page 183]{GrenanderSzego}, and also \cite[Chapter 6]{Varadhan}.
In the next section we derive an analogous formula for the variance of the optimal smoothing error when using both past and future values of $u_\ell$.

\section{Least-variance smoothing}\label{leastvariance}

Given the power spectrum $d\mu$ of a random process we seek the optimal linear smoothing filter using both past and future observations. It turns out that the variance of the smoothing error is the harmonic mean of the spectral density of the random process, i.e., it relates to the $0$th Fourier coefficient of the inverse of the spectral density of the process. This result will be used in the next section for the purpose of identifying the MR-spectra which are consistent with a finite set of autocorrelation samples.

\begin{thm}\label{prop:minvalue}{\sf
Let $d\mu$ be a bounded nonnegative measure on $[-\pi,\pi)$, let
$d\mu=d\mu_\s+fd\theta$ be the decomposition of $d\mu$ into its singular and absolutely continuous parts.
Then, the infimum of
\begin{equation}\label{minimalvalue}
\frac{1}{2\pi}\int_{-\pi}^{\pi}|\alpha(\theta)|^2d\mu(\theta)
\end{equation}
subject to the constraints
\begin{eqnarray}\label{constraint1}
\alpha(\theta)&\in& L_1,\\\frac{1}{2\pi}\int_{-\pi}^{\pi}\alpha(\theta)d\theta&=&1 \label{constraint2}
\end{eqnarray}
is equal to
\begin{equation}\label{infimum}
\left(\frac{1}{2\pi}\int_{-\pi}^{\pi} f(\theta)^{-1}d\theta\right)^{-1}
\end{equation}
when $f^{-1}\in L_1$, and zero otherwise.
}
\end{thm}

An important step in the proof of the theorem is provided by the following lemma.

\begin{lemma}\label{keylemma}\cite{GrenanderSzego} {\sf
Let $d\mu_\s$ be a bounded singular measure on $\theta\in I:=[-\pi,\pi)$ (i.e., the absolutely continuous part of $d\mu_\s$ is identically zero) and let $\epsilon_1$ be an arbitrary positive number. Then, it is always possible to decompose the interval $I$ into a finite number of intervals such that for a certain class $I_1$ of these intervals (i.e., their union) and for the complementary class $I_2=I\backslash I_1$, the following inequalities hold:
\begin{eqnarray*}
\frac{1}{2\pi}\int_{I_1}d\mu_\s(\theta)&<&\epsilon_1,\\
\frac{1}{2\pi}\int_{I_2}d\theta&<&\epsilon_1.
\end{eqnarray*}
}\end{lemma}

For a proof of Lemma \ref{keylemma} see \cite[page 7]{GrenanderSzego}. We now proceed with the proof of Theorem \ref{prop:minvalue}.

\begin{proof}{\em of Theorem \ref{prop:minvalue}:}
Assume first that $d\mu$ is absolutely continuous with no singular part. Given any positive number $\epsilon$ define
\[
\alpha_\epsilon(\theta):=\frac{\left( f(\theta)+\epsilon\right)^{-1}}{\frac{1}{2\pi}\int_{-\pi}^{\pi}\left( f(\theta)+\epsilon\right)^{-1}d\theta}.
\]
We note that $\alpha_\epsilon\in L_1$ (also in $L_2$ and in fact, it is even bounded and positive),
\[
\frac{1}{2\pi}\int_{-\pi}^{\pi}\alpha_\epsilon(\theta)d\theta=1,
\]
and we observe that
\begin{eqnarray*}
&&\frac{1}{2\pi}\int_{-\pi}^{\pi}|\alpha_\epsilon(\theta)|^2f(\theta)d\theta\\
&=&\frac{\frac{1}{2\pi}\int_{-\pi}^{\pi}\left( f(\theta)+\epsilon\right)^{-2}f(\theta)d\theta}{\left(\frac{1}{2\pi}\int_{-\pi}^{\pi}\left( f(\theta)+\epsilon\right)^{-1}d\theta\right)^2}\\
&< &\left(\frac{1}{2\pi}\int_{-\pi}^{\pi}\left( f(\theta)+\epsilon\right)^{-1}d\theta\right)^{-1}
\end{eqnarray*}
because $f(\theta)/(f(\theta)+\epsilon)<1$.
If $f^{-1}\not\in L_1$ then
\[
\lim_{\epsilon\to 0}\left(\frac{1}{2\pi}\int_{-\pi}^{\pi}\left( f(\theta)+\epsilon\right)^{-1}d\theta\right)^{-1}=0,
\]
whereas if $f^{-1}\in L_1$ the limit equals the expression given in (\ref{infimum}).
To prove our claims for the case where $d\mu$ is absolutely continuous, it remains to show that when $f^{-1}\in L_1$ the 
infimal value for (\ref{minimalvalue}) is never strictly less than (\ref{infimum}).

Continuing on, we assume that $f^{-1}\in L_1$.
We normalize $f^{-1}$ to have the identity as its $0$th Fourier coefficient
\[
\alpha_0=\frac{f^{-1}}{\frac{1}{2\pi}\int_{-\pi}^{\pi}f^{-1}(\theta)d\theta}
\]
and consider the perturbation
\[
\alpha=\alpha_0+\delta
\]
for an arbitrary $\delta\in L_1$ with vanishing $0$th Fourier coefficient (i.e., a $\delta\in L_1$ satisfying $\int_{-\pi}^\pi\delta(\theta)d\theta=0$ so that $\alpha$ satisfies (\ref{constraint2})). It readily follows that
\begin{eqnarray}\nonumber
&&\frac{1}{2\pi}\int_{-\pi}^{\pi}|\alpha(\theta)|^2d\mu(\theta)\\ \nonumber
&&\hspace*{-25pt} =\frac{1}{2\pi}\int_{-\pi}^{\pi}\left(|\alpha_0(\theta)|^2 + 2\delta(\theta)\alpha_0(\theta)+|\delta(\theta)|^2\right)f(\theta)d\theta\\
&&\hspace*{-25pt} =\frac{1}{2\pi}\int_{-\pi}^{\pi}|\alpha_0(\theta)|^2f(\theta)d\theta+
\frac{1}{2\pi}\int_{-\pi}^{\pi}|\delta(\theta)|^2f(\theta)d\theta, \label{laststep}
\end{eqnarray}
where for the last step we note that
\begin{eqnarray*}\int_{-\pi}^{\pi}\delta(\theta)\alpha_0(\theta)f(\theta)d\theta&=&
\int_{-\pi}^{\pi}\delta(\theta)\frac{f^{-1}(\theta)}{\frac{1}{2\pi}\int_{-\pi}^{\pi}f^{-1}(\theta)d\theta}f(\theta)d\theta\\
&=&\frac{1}{\frac{1}{2\pi}\int_{-\pi}^{\pi}f^{-1}(\theta)d\theta}\int_{-\pi}^{\pi}\delta(\theta)d\theta\\
&=&0.
\end{eqnarray*}
The first term in (\ref{laststep}) is precisely the claimed infimal value in (\ref{infimum}) and the second term is clearly nonnegative. This proves our claim in the case where $d\mu$ is absolutely continuous.

We now consider the case where
\[d\mu(\theta)=d\mu_\s(\theta)+f(\theta)d\theta\]
with $d\mu_\s$ a singular measure (always with respect to the Lebesgue measure). For an arbitrary $\epsilon>0$ we consider a decomposition of
\[
[-\pi,\,\pi)=I_1\cup I_2
\]
where 
\begin{eqnarray}\label{epsilon3a}
\int_{I_1}d\mu_\s(\theta)&<&\epsilon^3,\\\label{epsilon3b}
\int_{I_2}d\theta&<&\epsilon^3.
\end{eqnarray}
That such a decomposition exists follows from Lemma \ref{keylemma} taking $\epsilon_1=\epsilon^3$ in the statement of the lemma. Now let $\chi_{I_1}$ denote the characteristic function of $I_1$ which takes the value $1$ when $\theta\in I_1$ and zero otherwise, and set
\begin{equation}\label{eq:star}
\alpha_\epsilon=\frac{(f+\epsilon)^{-1}\chi_{I_1}}{\frac{1}{2\pi}\int_{-\pi}^{\pi}(f(\theta)+\epsilon)^{-1}\chi_{I_1}(\theta)d\theta}
\end{equation}
which is in $L_1$ and has the identity as its $0$th Fourier coefficient. Then
\begin{eqnarray*}
\frac{1}{2\pi}\int_{-\pi}^{\pi}|\alpha_\epsilon(\theta)|^2d\mu(\theta)&=&
\frac{1}{2\pi}\int_{-\pi}^{\pi}|\alpha_\epsilon(\theta)|^2d\mu_\s(\theta)\\
&&\hspace*{-3cm}+\frac{1}{2\pi}\int_{-\pi}^{\pi}   
\frac{(f+\epsilon)^{-2}\chi_{I_1}}{\left(\frac{1}{2\pi}\int_{-\pi}^{\pi}(f(\theta)+\epsilon)^{-1}\chi_{I_1}(\theta)d\theta\right)^2}
f(\theta)d\theta.
\end{eqnarray*}
The first term on the right hand side is bounded above by
\[
\frac{1/\epsilon^2}{\left(\frac{1}{2\pi}\int_{-\pi}^{\pi}(f(\theta)+\epsilon)^{-1}\chi_{I_1}(\theta)d\theta\right)^2}\epsilon^3
\]
which decays to $0$ with $\epsilon$, whereas the second term is bounded above by
\[
\frac{1}{\left(\frac{1}{2\pi}\int_{-\pi}^{\pi}(f(\theta)+\epsilon)^{-1}\chi_{I_1}(\theta)d\theta\right)}
\]
which in the limit recovers the claimed bound (\ref{infimum}). The earlier argument for the case of absolutely continuous $d\mu$ applies and shows that this bound is in fact the correct value for the infimum and that no lower value is possible.
\end{proof}

\begin{remark}{It is clear from the proof that if $f^{-1}\in L_1$ and $d\mu=fd\theta$ is absolutely continuous, then 
\[
\alpha_0=\frac{f^{-1}}{\frac{1}{2\pi}\int_{-\pi}^{\pi}f^{-1}(\theta)d\theta}
\]
is the unique optimal solution which achieves the minimal value
\[
\left(\frac{1}{2\pi}\int_{-\pi}^{\pi} f(\theta)^{-1}d\theta\right)^{-1}
\]
for 
\[
\|\alpha\|_{d\mu}^2:=\frac{1}{2\pi}\int_{-\pi}^{\pi}|\alpha(\theta)|^2d\mu(\theta)
\]
subject to $\alpha\in L_1$ with $0$th Fourier coefficient the identity.
Thus, if
\[
\alpha_0(\theta)\sim \ldots +\rho_{-1}e^{-j\theta}+1+\rho_1e^{j\theta}+\ldots
\]
and $\rho_\ell$, $\ell=\pm 1,\pm 2,\ldots$ the corresponding Fourier coefficients, then
\[
\hat{u}_0=-\sum_{\ell\neq 0} \rho_\ell u_{-\ell}
\]
is the optimal in the least variance sense estimate for $u_0$, and $u_\ell$ is a random process with
$d\mu$ as its power spectrum. In this case the infimum is achieved, and hence it represents the minimum variance of the error. When the power spectrum has either singular part or $f^{-1}\not\in L_1$, then $\alpha_\epsilon$ as in (\ref{eq:star}) provides suboptimal solutions. This is completely analogous to the Szeg\"{o}-Kolmogorov setting where optimal one-step ahead predictors (which use only past observations) exist when $\log f\in L_1$ otherwise the least variance is not attained but can be gotten arbitrarily closely \cite[Chapter II]{Geronimus}.}\end{remark}

\begin{remark} It is interesting to observe that while the minimal variance of a smoothing error for a random process having $f$ as spectral density is the {\em harmonic mean}
\[m_{-1,f}:=\left(\frac{1}{2\pi}\int_{-\pi}^{\pi} f(\theta)^{-1}d\theta\right)^{-1}\]
of the values of $f$ on the $[-\pi,\,\pi)$, the minimal variance of the optimal one-step-ahead predictor using only past observations is the {\em geometric mean} (see \cite[page 183]{GrenanderSzego}, \cite[Chapter 6]{Varadhan})
\[ m_{0,f}:=\exp \left(\frac{1}{2\pi}\int_{-\pi}^{\pi}\log\left(f(\theta)\right)d\theta\right).\]
The former is the inverse of $\bJ (d\mu/d\theta)$
whereas the latter is exponential of $-\bI (d\mu/d\theta)$.
Naturally, $m_{-1,f}\leq m_{0,f}$ (see also \cite[page 23]{BB}). This ordering is clear
from the interpretation of the two quantities as variances of best predictors which use ``past$+$future'' and ``only past'' observations, respectively.
\end{remark}

\section{On deterministic processes: an example}

It may be rather surprising, at first glance, that the value of a random process with power spectral density
\begin{equation}\label{fo}
f_o(\theta)=|1-e^{j\theta}|^2= 2-2\cos(\theta)
\end{equation}
can be predicted at any given point with arbitrarily small variance, when both past and future observations are available. Yet this is the case, and this is due to the fact that $f_o^{-1}\not\in L_1$ (equivalently $m_{-1,f_o}=0$). This example highlights the difference between ``deterministic processes'' in the sense of $m_{-1,f}=0$ and those in the sense of Szeg\"{o}-Kolmogorov which are characterized by $m_{0,f}=0$ or, equivalently, by $\log f\not\in L_1$ instead.

For our particular example, the fact that $f_o^{-1}\not\in L_1$ follows from the divergence of
\[
\int_\epsilon^\pi \frac{1}{1-\cos(\theta)}d\theta >\int_\epsilon^\pi \frac{1}{\theta^2}d\theta
\]
as $\epsilon\to 0$. On the other hand, the fact that $\log f_o\in L_1$ can be seen as follows. Since
$g(z):=1-z$ is analytic and does not vanish in $\mD$, $\log |g|$ is harmonic and
\[
\frac{1}{2\pi}\int_{-\pi}^\pi \log(|g(re^{j\theta})|)d\theta =\log(|g(0)|)=0
\]
for any value of $r\in[0,1)$. Therefore, the integral of the logarithm of $\lim_{r\to 1} |g(re^{j\theta})|$ also vanishes, and the same applies to $f_o(\theta)=\lim_{r\to 1} |g(re^{j\theta})|^2$.

In the rest of this section we explain how a random process corresponding to $f_o$ can be predicted with vanishingly small variance from the combined past and future record. We do so, for didactic purposes, by sketching a specialized and more direct construction than that of Section \ref{leastvariance}.

Consider a realization $\{u_k,\,k\in\mZ\}$ of a random process corresponding to $f_o$ as follows:
\[
u_k=w_k-w_{k-1}
\]
where $\{w_k,\,k\in\mZ\}$ is  a sequence of independent, identically distributed, random variables with zero mean and unit variance (i.e., a white-noise process). We assume that ``past'' ($\{u_k,\,k<0\}$) as well as ``future'' ($\{u_k,\,k>0\}$) observations are available, and that we wish to estimate the ``present'' $u_0=w_0+w_{-1}$ based on this two-sided observation record. Then,
\newcommand{\bu}{{\bf u}}
\newcommand{\bw}{{\bf w}}
\[
\bu_{_{<0}}:=\left[\begin{array}{c}u_{-1}\\u_{-2}\\u_{-3}\\\vdots\end{array}\right]
=\left[\begin{array}{cccc }1 & -1 & 0 & \ldots\\
                                        0 & 1 & -1 & \ddots\\
                                        0 & 0 & 1 &\ddots  \\
                                        \vdots & \vdots &\ddots &\ddots \end{array}\right]
\left[\begin{array}{c}w_{-1}\\w_{-2}\\w_{-3}\\\vdots\end{array}\right] 
\]
and
\[
\bu_{_{>0}}:=\left[\begin{array}{c}u_{1}\\u_{2}\\u_{3}\\\vdots\end{array}\right]
=\left[\begin{array}{cccc }-1 & 1 & 0 & \ldots\\
                                        0 & -1 & 1 & \ddots\\
                                        0 & 0 & -1 &\ddots  \\
                                        \vdots & \vdots &\ddots &\ddots \end{array}\right]
\left[\begin{array}{c}w_{0}\\w_{1}\\w_{2}\\\vdots\end{array}\right] 
\]
In both cases the mapping is Toeplitz, and identical except for a sign change.
Let now
\[
v:=\left[\begin{array}{cccc}1&(1-\epsilon) & (1-\epsilon^2) &\ldots\end{array}\right],
\]
and for $1>\epsilon>0$ and define
\begin{eqnarray*}
\hat{w}_{-1} &:=& v\bu_{_{<0}} =w_{-1} +\epsilon(1-\epsilon)\sum_{k=-2}^{-\infty} \epsilon^{-k+2} w_k\\
\hat{w}_{0} &:=& -v\bu_{_{>0}} =w_{0} +\epsilon(1-\epsilon)\sum_{k=1}^\infty \epsilon^{k-1} w_k\\
\hat{u}_0&:=&\hat{w}_{0}-\hat{w}_{-1}.
\end{eqnarray*}
Each of the above can be taken as an estimator for the corresponding un-hatted variable.
The variance of estimation in all cases can be made arbitrarily small with appropriately small choice for $\epsilon$. This justifies our claim.

\section{Proofs of Theorems \ref{thm2} and \ref{computation}}\label{proofs}

Due to the strict concavity of the inversion map $x\mapsto 1/x$ on $\mR_+$, $\bJ(\cdot)$ is also a strictly concave functional on (non-negative) density functions. We first show that a spectral density $f_\me$ of the form claimed in Theorem \ref{thm2} is indeed a minimizer of $\bJ(\cdot)$ subject to the moment constraints
\begin{equation}\label{fmoment}
R_k=\frac{1}{2\pi}\int_{-\pi}^{\pi}e^{-jk\theta}f(\theta)d\theta,\mbox{ for }k=0,1,\ldots,n.
\end{equation}
Existence of suitable values for the corresponding parameters requires proving Theorem \ref{computation} next,  which claims that these values correspond to an attractive equilibrium of a certain differential equation. 
The form of $f_\me$ ensures stationarity and hence, due to the strict concavity of $\bJ(\cdot)$, it ensures that this is indeed the unique extremal point. Finally, we revisit the optimizaton problem and consider measures with possible singular part. The singular part does not affect the value of $\bJ(\cdot)$, but the fact that a singular part is allowed, relaxes the constraint (\ref{fmoment}) to (\ref{moment}). Yet, as we will see, $f_\me$ is still the minimizer and, hence, the extremal spectral measure $d\mu$ cannot have a singular part. In the end, we return to the remaining claims in Theorem \ref{thm2} regarding properties of the minimizer.

\subsection{Functional form of minimizer}\label{functionalform}
Consider first the problem of minimizing $\bJ(f)$ with $f$ constrained to satisfy (\ref{fmoment}). If
\[
\lambda:=\left[\begin{array}{ccccc}\lambda_{-n}&\ldots&\lambda_0&\ldots&\lambda_n\end{array}\right]
\]
denotes a vector of Lagrange multipliers, the corresponding Lagrangian is
\newcommand{\cL}{{\mathcal L}}
\begin{eqnarray*}
\cL(f,\lambda)&:=&\bJ(f)-\lambda(\rR_1-\frac{1}{2\pi}\int_{-\pi}^\pi G(e^{j\theta}) f(\theta)d\theta)
\end{eqnarray*}
where $\rR_1,G$ are defined in the statement of Theorem \ref{computation}.
If we set the variation
\[
\delta\cL(f,\lambda;\delta f)= \frac{1}{2\pi}\int_{-\pi}^\pi \left(\frac{-1}{f(\theta)^2}+\lambda G(e^{j\theta})\right) \delta f(\theta)d\theta
\]
identically equal to zero for all perturbations $\delta f$ (assuming that $f>0$ and hence $\delta f$ unconstrained), then we conclude that 
\begin{equation}\label{eq:sufficient}
f(\theta)=\frac{1}{\sqrt{\lambda G(e^{j\theta})}},
\end{equation}
which is the form claimed in Theorem \ref{thm2} for $f_\me$. Our next step is to prove that, provided $\bR_n>0$, there always exists such a density function which satisfies (\ref{fmoment}) and that the trigonometric polynomial $\lambda G(e^{j\theta})$ is in fact strictly positive.

\subsection{Proof of Theorem \ref{computation}}\label{proofofcomputation}
We follow the formalism in \cite{IT} for solving moment problems.
We denote by $\fR$ the positive cone
\[
\fR:=\{\rR\;:\;\rR=\intpi G(e^{j\theta}) d\mu(\theta),\mbox{ where }d\mu\geq 0\}
\]
and by $\fK$ the dual cone
\[
\fK:=\{\lambda \;:\;\lambda G(e^{j\theta}) \geq 0\mbox{ for }\theta\in[-\pi,\pi]\}.
\]
Both are subsets of $\mR\times \mC^{2n}$ since their ``$0$th'' entries $R_0,\lambda_0\in\mR_+$ while the remaining entries $R_\ell,\lambda_\ell\in\mC$ ($\ell=\pm 1,\pm 2,\ldots,\pm n$). Also, both are convex. The interior of $\fR$ is denoted by ${\rm int}(\fR)$ and the interior of the dual cone, which consists of all vectors $\lambda$ such that the trigonometric polynomial $\lambda G$ is strictly positive on the unit circle, is denoted by $\fK_+$. The Jacobian $\frac{\partial H}{\partial \lambda}$ of the mapping
\[
H \;:\; \fK_+\to {\rm int}(\fR) \;:\; \lambda\mapsto \intpi G(e^{j\theta})\frac{1}{\sqrt{\lambda G(e^{j\theta})}} \, d\theta
\]
between Lagrange vectors and moments is given in (\ref{Jacobian}) and is denoted by $M(\lambda)$. As long as $\lambda\in\fK_+$ the Jacobian is an invertible matrix. Our goal is to find a value for $\lambda$ so that condition (ii) of Theorem \ref{computation} holds. We do this as follows.

We begin with $\lambda_0$ as in Theorem \ref{computation} for which we readily observe that $\lambda_0 G\equiv 1>0$. It follows that
\[
\rR_0:=\intpi G(e^{j\theta})\frac{1}{\sqrt{\lambda_0 G(e^{j\theta})}}\,d\theta \in{\rm int}(\fR).
\]
Since $\bR_n>0$, we also know that $\rR_1\in{\rm int}(\fR)$. Since ${\rm int}(\fR)$ is convex and $\rR_1,\rR_0\in{\rm int}(\fR)$, the
interval $[\rR_0,\,\rR_1]\subset {\rm int}(\fR)$, i.e.,
\begin{equation}\label{interval}
\rR_\tau:=\tau \rR_1 +(1-\tau)\rR_0
\end{equation}
belongs to ${\rm int}(\fR)$ for all $\tau\in[0,1]$. The key idea is now to trace $\rR_\tau$ by following corresponding values for $\lambda_\tau$ in the dual cone. This is not always possible. It depends on the functional form for the sought spectral density function $f$. The critical issue that may prevent such path-following in the dual space is whether any $\lambda$ in the boundary of $\fK_+$ maps onto a point in the interior of $\fR$. When this happens, there are interior points in $\fR$ which do not admit the assumed representation. We will see below that this does not happen for the functional form $1/\sqrt{\lambda G}$ and hence, that the plan we have outlined applies. We discuss these key steps/facts next.

The moments $\rR_\tau$, $\tau\in[0,1]$, satisfy the differential equation
\begin{equation}\label{diffeq0}
\frac{d\rR_\tau}{d\tau}=\rR_1-\rR_0
\end{equation}
as follows readily from (\ref{interval}).
Then the dual parameters $\lambda(\tau)$ satisfy
\begin{equation}\label{diffeq}
\frac{d\lambda(\tau)}{d\tau}=M(\lambda)^{-1}(\rR_1-\rR_0),
\end{equation}
as long as $\lambda(\tau)$ remains in the interior of $\fK_+$ ---in which case $M(\lambda)$ is invertible being the (inverse of the) autocorrelation matrix of a positive spectral density function. We claim that this is always the case. To prove it, assume that the contrary is true and that $[0,\tau_0)$ is a maximal subinterval of $[0,1]$ for which $\lambda(\tau)\in\fK_+$ for $0\leq \tau<\tau_0$. Thus, the family of positive trigonometric polynomials
\[
\{\lambda(\tau)G(e^{j\theta})\;:\;\tau\in[0,\tau_0)\}
\]
has either a limit point on the boundary of $\fK_+$ or it grows unbounded. In either case we will draw a contradiction.

In the first case, there must exist an accumulation point $\hat{\lambda}$
for which $\hat{\lambda}G(e^{j\theta})$ vanishes
on the unit circle. But then $\hat{\lambda}G(e^{j\theta})$, which is a nonnegative trigonometric polynomial, must have a double root at some point $e^{j\theta_0}$. Therefore
\begin{equation}\label{nonintegrable}
\frac{1}{\sqrt{\hat{\lambda}G(e^{j\theta})}},
\end{equation}
which has at least a single pole at $e^{j\theta_0}$, is not integrable. The assertion that the inverse of the square root of a nonnegative trigonometric polynomial which vanishes on the circle is not integrable is elementary. It suffices to consider a typical case, such as $1-\cos(\theta)$, where $\frac{1}{\sqrt{1-\cos(\theta)}}=\frac{1}{\sqrt{2}|\sin(\theta/2)|}>\frac{\sqrt{2}}{|\theta|}$ is clearly not integrable---the general case is similar. The nonintegrability of (\ref{nonintegrable}) implies that
the family of vectors
\[
\{\intpi G(e^{j\theta})\frac{1}{\sqrt{\lambda(\tau)G(e^{j\theta})}}\,d\theta \;:\;0\leq\tau<\tau_0\}
\]
is unbounded, in contradiction to the assumption that the image of 
$\{\lambda(\tau)\;:\;0\leq \tau<\tau_0\}$ under $H$ is the subset
\[\{\rR_\tau\;:\;0\leq \tau<\tau_0\}
\]
of the bounded interval $[\rR_0,\rR_1]$.

We now draw a contradiction for the second case.
We assume that $\lambda(\tau)$ grows unbounded as $\tau\to\tau_0$.
It follows that there is sequence $\tau_i\in[0,\tau_0)$, $i=1,2,\ldots$ such that $\tau_i\to \tau_0$ and
$\|\lambda(\tau_i)\|\to \infty$ while the unit-length vectors
\[
\hat{\lambda}_i:=\frac{\lambda(\tau_i)}{\|\lambda(\tau_i)\|}\to \hat{\lambda}\in\fK
\]
converge as $i\to\infty$, with $\|\cdot\|$ being the Euclidean norm.
At the same time, the sequence $\rR_{\tau_i}=H(\lambda(\tau_i))$, $i=1,2,\ldots$, converges to
$R_{\tau_0}\in{\rm int}(\fR)$.
But any interior point $\rR\in\fR$ is characterized by the property that
the functional
\[ \fC_{\rR}\;:\;\fK\to\mR_+\;:\;\lambda\mapsto \lambda \rR\]
is strictly positive (e.g., see \cite[Proposition 3]{IT}). (This is due to the fact any such $\rR$ assumes a representation
$\intpi G(e^{j\theta})f(\theta)d\theta$ for some strictly positive density function $f(\theta)$.) 
On the other hand, returning to the sequence $R_{\tau_i}$ $i=1,2,\ldots$, we observe that
\begin{eqnarray*}
\fC_{\rR_{\tau_i}}\;:\; \hat{\lambda}_i \mapsto \hat{\lambda}_iR_{\tau_i}&=&\intpi \frac{\hat{\lambda}_iG(e^{j\theta})}{{\sqrt{\lambda(\tau_i) G(e^{j\theta})}}}\,d\theta\\
&=&\intpi \sqrt\frac{{\hat{\lambda}_i G(e^{j\theta})}}{\|\lambda(\tau_i)\|}\,d\theta
\end{eqnarray*}
tends to $0$ as $\|\lambda(\tau_i)\|$ grows unbounded. Therefore, the functionals $\fC_{\rR_{\tau_i}}$, $i\to\infty$, are not uniformly bounded away from zero. Yet, their limit $\fC_{\rR_{\tau_0}}$ is, due to the fact that $\rR_{\tau_0}\in{\rm int}(\fR)$. This is a contradiction. Therefore (\ref{diffeq}) can be integrated over the complete interval $[0,1]$
and $\lambda(\tau)$ remains bounded and in the interior of the dual cone (i.e., the trajectory lies in $\fK_+$). We identify $\lambda(1)=\lambda_\mr$.

We now re-scale the independent variable in (\ref{diffeq0}-\ref{diffeq}) by replacing $\tau$ with $t=-\log(1-\tau)$.
We simplify notation and denote $\rR_{\tau(t)}$ by $\rR_t$ and $\lambda(t(\tau))$ by $\lambda(t)$. 
Using $\frac{\partial \tau}{\partial t}=1-\tau$ and $\rR_1-\rR_0=\frac{1}{1-\tau}(\rR_1-\rR_\tau)$, 
we rewrite (\ref{diffeq0}) as
\[
\frac{d\rR_t}{dt}=\rR_1-\rR_t, \mbox{ for }t\in[0,\infty),
\]
and (\ref{diffeq}) as
\begin{equation}\label{diffeq2}
\frac{d\lambda(t)}{dt}=M(\lambda(t))^{-1}(\rR_1-\rR_t),
\end{equation}
where, as usual, $\rR_t=\intpi G\frac{1}{\lambda(t)G}\,d\theta$. We have now established claims (i) and (ii) of Theorem \ref{computation}. I.e., we have shown that as $t\to\infty$ in (\ref{diffeq2}) the trajectory $\lambda(t)$ converges in $\fK_+$, and that the limit point $\lambda_\mr$ is such that (\ref{moment}) holds. Claim (iii) of the theorem follows from the concavity of $\bJ(\cdot)$. More specifically, the functional form of $f_\mr$ guarantees that it is a minimizer of $\bJ(\cdot)$. There can only be one such minimizer since $\bJ(\cdot)$ is strictly concave.

\subsection{Proof of Theorem \ref{thm2}}\label{proofofthm2}

Define first the column vector
\[
g(e^{j\theta}):=\left[\begin{array}{cccc}1&e^{-j\theta}&\ldots&e^{-jn\theta}\end{array}\right]^\prime.
\]
Assuming that $d\mu=d\mu_\s+fd\theta$ with $d\mu_s$ a singular measure and $fd\theta$ the absolutely continuous part of $d\mu$, the minimization of $\bJ(f)$ subject to (\ref{moment}) is equivalent to minimization of $\bJ(f)$ subject to
\begin{equation}\label{inequalityconstraint}
\bR_n\geq \intpi g(e^{j\theta})f(\theta)g(e^{j\theta})^*d\theta.
\end{equation}
The corresponding Lagrangian is now
\begin{eqnarray}
\cL_o(f,\Lambda)&:=&\bJ(f(\theta))+ \label{Lo}\\
&& \hspace*{-2cm}+{\rm trace}\left( \Lambda \left(\bR_n- \intpi g(e^{j\theta})f(\theta)g(e^{j\theta})^*d\theta\right) \right)\nonumber\\
&=&\bJ(f(\theta))+{\rm trace}(\Lambda\bR_n) \nonumber\\
&& \label{Lo2}
-\intpi\left(g(e^{j\theta})^*\Lambda g(e^{j\theta})\right)f(\theta)d\theta
\end{eqnarray}
The Lagrange multiplier $\Lambda$ is a matrix which has a Toeplitz structure. (To see this note that any possible component of $\Lambda$ which is orthogonal to the subspace of Toeplitz matrices has no effect since it vanishes when taking the inner product ${\rm trace}(\Lambda T)$ for any Toeplitz matrix $T$ as done in (\ref{Lo}).) The minimizer $f$ would correspond to a measure $d\mu$ with a nontrivial singular part only if the equality constraint in (\ref{inequalityconstraint}) is not active. For this to be the case, the multiplier
\[
g(e^{j\theta})^*\Lambda g(e^{j\theta})
\]
of $f(\theta)$ in (\ref{Lo2}) must vanish at least for some values of $\theta$.
However, the correspondence
\[
\Lambda=\left[\begin{array}{cccc}
\frac{1}{n+1}\lambda_0 &\frac{1}{n}\lambda_1&\ldots&\frac{1}{1}\lambda_n\\
\frac{1}{n}\lambda_1&\frac{1}{n+1}\lambda_0 &\ldots&\frac{1}{2}\lambda_{n-1}\\
\vdots&\vdots&\ddots&\vdots\\
\frac{1}{1}\lambda_{-n}&\frac{1}{2}\lambda_{-n+1}&\ldots&\frac{1}{n+1}\lambda_0
\end{array}\right]
\]
shows that in fact $\cL_o(f,\Lambda)=\cL(f,\lambda)$, i.e., it is the same Lagrangian as in Section \ref{functionalform}. The value for the Lagrange multipliers in the latter, as identified in Section \ref{proofofcomputation}, are such that $\lambda_\mr G(e^{j\theta})$ is a positive trigonometric polynomial. This polynomial is precisely the multiplier of $f(\theta)$ in (\ref{Lo2}) and is strictly positive for all $\theta\in[-\pi,\pi]$. Hence, the equality constraint in (\ref{inequalityconstraint}) is active for the extremal $f$ of the relaxed problem corresponding to (\ref{Lo}).
Then, the analysis in Section \ref{functionalform} applies. Therefore, the minimizer corresponds to an absolutely continuous power spectral distribution $d\mu_\mr=f_\mr(\theta)d\theta$ which is of the form claimed in the theorem.

We now address the remaining claims in the theorem regarding the variance of the smoothing error for the corresponding random process. 
Given the expression for $f_\mr$ which is the square root of the inverse of a positive trigonometric polynomial, the form of the optimal smoothing filter for the corresponding random process is provided by Theorem \ref{prop:minvalue}. It is a consequence of the same theorem that
the variance of the optimal smoothing error $\cE_{d\mu_\mr}\{|u_0 -\hat{u}_{0|{\rm past \;\&\; future}}|^2\}$ is
precisely the inverse of the $\bJ$-functional evaluated at $f_\mr$, i.e.,
\[
\left(\bJ(f_\mr)\right)^{-1}.
\]
The last part of the theorem is also immediate since
\[
\min_{\beta_k,\;k\neq 0}\left\{ \cE_{d\mu}\{|u_0 -\sum_{k\neq 0} \beta_k u_{-k}|^2\}\;:\; \mbox{ } \mbox{(\ref{moment}) holds}\right\}
\]
is $(\bJ(d\mu/d\theta))^{-1}$ for any spectral measure consistent with (\ref{moment}). But $d\mu_\mr$ is the unique maximizer of this inverse.

\section{On spectral analysis: an example}
For illustration purposes, we compare the power spectra $f_\me$ and $f_\mr$ given in Theorems \ref{MEtheorem} and \ref{thm2} for a basic example.
We begin by evaluating the first $4$ autocorrelation moments for the following spectral density:
\[
f_{\rm true}(\theta)= 1+\frac{2}{5}\cos(\theta)+\delta(\theta-\frac{1}{2})+\delta(\theta+\frac{1}{2}).
\]
Here, for convenience, we depart slightly from our earlier notation and incorporate the singular part of the power spectrum into the ``spectral density'' as a sum of two Dirac functions ---the distributions $\delta(\theta-\theta_0)$ for $\theta_0=\pm\frac{1}{2}$. Thus, the absolutely continuous part of the power spectrum is made up of only the continuous portion $\left(1+\frac{2}{5}\cos(\theta)\right)d\theta$ of $f_{\rm true}(\theta)d\theta$.
The corresponding random process consists of a random moving average component generated by
\[
u^{MA}_k=w_k+\frac{1}{2}w_{k-1}
\]
with $w_k$ a white-noise process with variance $1/(1+1/4)$ (normalized so that $\cE\{|u_k^{MA}|^2\}=1$), and a deterministic sinusoidal component at frequency $\theta_{o}=1/2$ [rad/unit of time]. The first $4$ samples of the autocorrelation function of
\[
u_k=u_k^{MA}+2\sin(\frac{k}{2}+\phi)
\]
(with $\phi$, say, uniformly distributed on $[-\pi,\,\pi]$) can be readily computed and are as follows:
\begin{eqnarray*}
&&\left[\begin{array}{cccc}R_0 &R_{\pm 1}&R_{\pm 2}&R_{\pm 3}\end{array}\right] \\
&&=\left[\begin{array}{cccc}3.0000   & 2.1552  &  1.0806 &   0.1415\end{array}\right].
\end{eqnarray*}
The corresponding Toeplitz matrix ${\bf R}_3$ is positive definite, and as a result, there is a nontrivial family of power spectra which are consistent with the autocorrelation data --$d\mu(\theta)=f_{\rm true}(\theta)d\theta$ is only of them.

\begin{figure}[htb]\begin{center}
\includegraphics[totalheight=7cm]{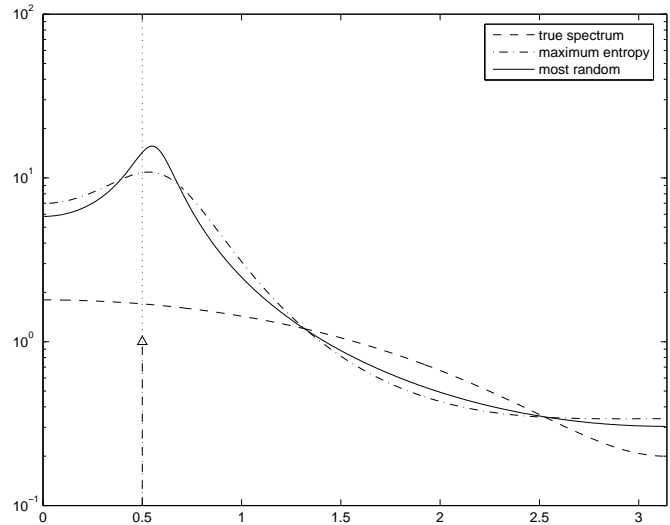}\\
\caption{Power spectra consistent with $R_0,\,R_{\pm 1},\,R_{\pm 2},\,R_{\pm 3}$.}\label{powerspectra}\end{center} \end{figure}

Figure \ref{powerspectra} shows the three particular power spectra that concern us here. First, the ``moving-average $+$ sinusoids'' power spectrum described above is shown with a dashed line ($- - -$). Then, a ME-power spectrum which is consistent with ${\bf R}_3$ and obtained following the maximum entropy ansatz is shown with a dash-dotted line ($- \cdot -$). Finally, the MR-power spectrum corresponding to the least smooth process is shown with a continuous line (-------).

{\em All three power spectra shown are consistent with the covariance data}. Hence, there is no suggestion that one should be preferable. They all describe the same data. A selection can only be based on either prior information or a prejudice ---this is where an ``ansatz'' becomes relevant. Had we known that the ``true'' spectrum originates from a moving average component plus a minimal number of sinusoids, we could have recovered the exact power spectrum from the covariance data following e.g., \cite{CAR}. Of course, such knowledge is rarely available and one is called to use other insights. Thence, if the power spectrum and a model for the process is to be used for prediction purposes, the maximum entropy option is quite natural since it represents the relevant ``worst-case senario.'' However, if the model is to be used for filling in gaps in records, then the MR-option is the appropriate ``worst-case senario.'' Then, if our goal is to simply identify features in the power spectrum, either may be appropriate.

Using (\ref{orthogonalpoly}) we determine that
\begin{eqnarray}\label{MaxEnt}
f_\me(\theta)&=&
\frac{k_\me^2}{|1+a_1 e^{j\theta} +a_2 e^{2j\theta} +a_3e^{3j\theta}|^2}
\end{eqnarray}
with $k_\me=1.2732$, and
\begin{eqnarray*}
&& \left[\begin{array}{ccc}a_1 &a_2 &a_3\end{array}\right]\\
&& =\left[\begin{array}{ccc}-0.9026   & 0.1829 &   0.1465\end{array}\right].
\end{eqnarray*}
On the other hand, following Theorem \ref{computation} we compute
\begin{eqnarray}\label{mostrandom}
f_\mr(\theta)&=&
\frac{1}{
\sqrt{\sum_{\ell=-3}^{3}\lambda_\ell e^{\ell j\theta}}}\\\nonumber
&=&\frac{\kappa^2}{\sqrt{|1+{\hat{a}}_1 e^{j\theta} +{\hat{a}}_2 e^{2j\theta} +{\hat{a}}_3e^{3j\theta}|^2}}\\ \label{mostrand2}
&=&\frac{\kappa^2}{|1+{\hat{a}}_1 e^{j\theta} +{\hat{a}}_2 e^{2j\theta} +{\hat{a}}_3e^{3j\theta}|}
\end{eqnarray}
with
\begin{eqnarray*}
&&\left[\begin{array}{cccc}\lambda_0 &\lambda_{\pm 1}&\lambda_{\pm 2}&\lambda_{\pm 3}\end{array}\right]\\
&&=\left[\begin{array}{cccc} 3.4942  & -2.5690  &  0.9598 &  -0.1231\end{array}\right],
\end{eqnarray*}
or, equivalently,  $\kappa =1.2732$ and
\begin{eqnarray*}
&&\left[\begin{array}{ccc}{\hat{a}}_1 &{\hat{a}}_2 &{\hat{a}}_3\end{array}\right]\\
&&=\left[\begin{array}{ccc}   -1.7673  &  1.1795  & -0.1956\end{array}\right].
\end{eqnarray*}
Here, again, we depart slightly from our earlier notation so as to compare the coefficients more directly to the ME-spectral density. The parameters $b_\ell$ and $\rho_\ell$ as in Theorem \ref{thm2} for $f_\mr$ and smoothing filter, respectively, can be readily determined from the above.

\begin{figure}[htb]\begin{center}
\includegraphics[totalheight=7cm]{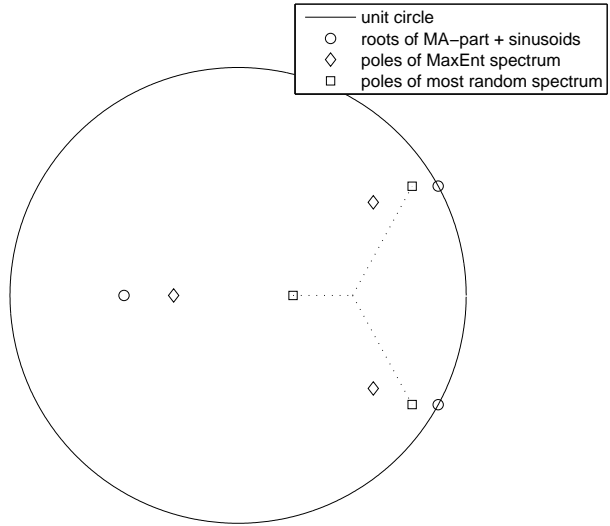}
\caption{Poles/zeros of $f_{\rm true}$, and singularities of $f_\me$, and $f_\mr$.}\label{fig:poles}\end{center} \end{figure}

Figure \ref{fig:poles} marks the zero of the moving average component of $f_{\rm true}$ (inside $\mD$) along with the location of the two spectral lines (on the unit circle) with ``o''. The poles of the ME-spectrum are marked with a ``$\diamond$'' and the fractional poles of the MR-spectrum with a ``{\footnotesize $\Box$}''.

\begin{figure}[htb]\begin{center}
\includegraphics[totalheight=7cm]{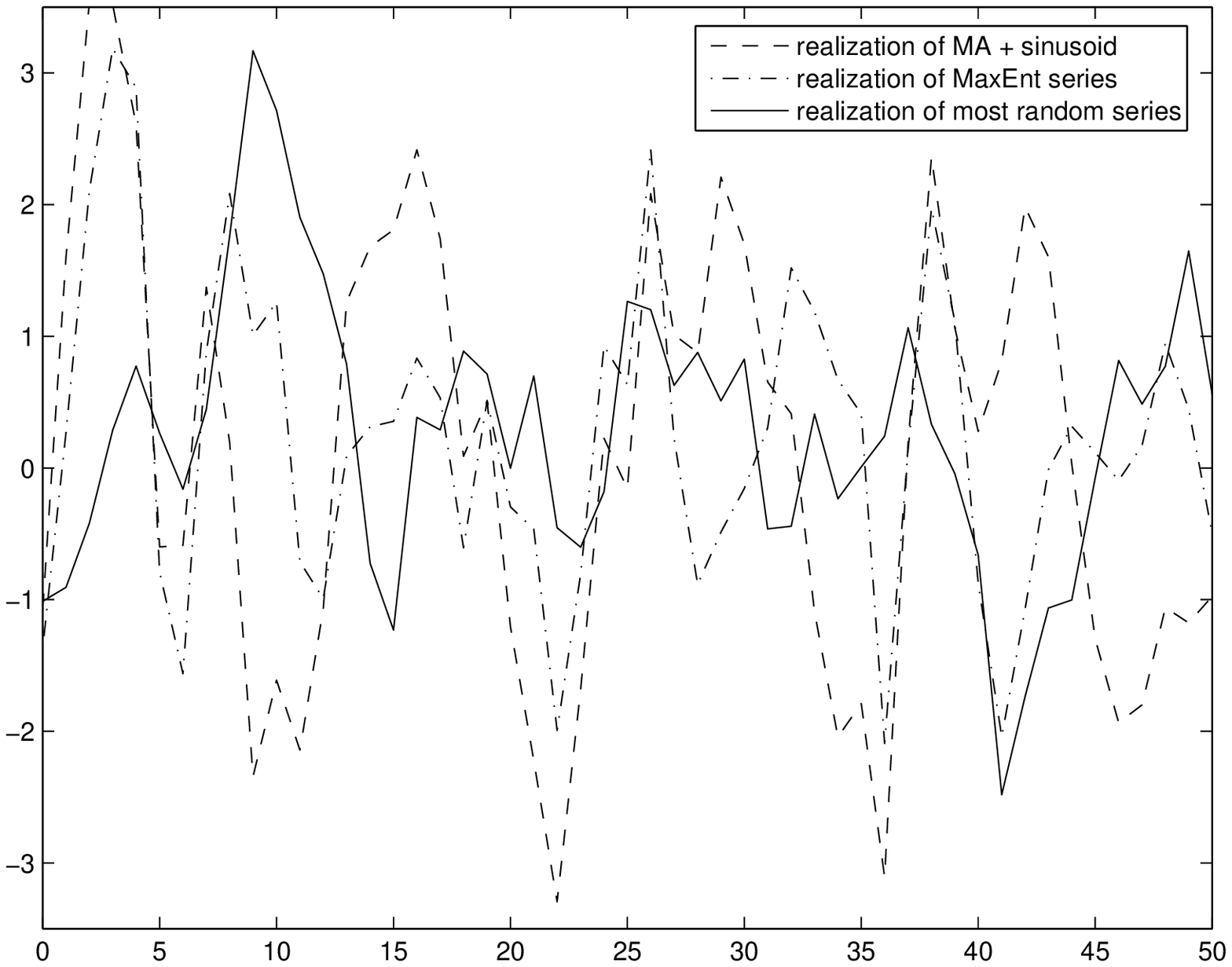}
\caption{Realizations of time-series according to $f_{\rm true}$, $f_\me$, and $f_\mr$.}\label{fig:simulation}\end{center} \end{figure}

Figure \ref{fig:simulation} presents realizations of time-series corresponding to $f_{\rm true}$, $f_\me$, and $f_\mr$. The one corresponding to $f_{\rm true}$ is generated by a Markovian moving-average model plus a sinusoidal component with a random phase. The time-series corresponding to $f_\me$ is generated by a Markovian autoregressive model as usual. Finally, the time-series corresponding to $f_\mr$ is generated by a suitable discretization of the standard spectral representation (stochastic integral)
\[
u_\ell=\frac{1}{2\pi}\int_{-\pi}^\pi e^{j\ell\theta}dv(\theta)
\]
where $dv(\theta)$ is a zero-mean white noise process for $\theta\in[-\pi,\,\pi)$ such that $\frac{d}{d\theta}\cE\{|v(\theta)|^2\}=f_\mr(\theta)$, see e.g., \cite[page 183]{GrenanderSzego}. There is not apparent observational feature distinguishing these three realizations, at least over the window where they have been drawn, and hence, they are produced here only to satisfy curiosity.

\section{Concluding remarks}

The present study sought to explore the issue of the time-arrow in the context of the maximum entropy ansatz. When the index of a random process designates a variable other than time, the principle can be called into question. A more abstract version of seeking spectra {\em maximally noncommittal to unavailable data}, such as gaps in a record, suggests other alternatives, including the one studied herein.

At the moment, the information theoretic significance of $\bJ(d\mu(\theta)/d\theta)$ is still under consideration. However, it is clear that, in the same way that entropy rates relate to a level of ``surprise'' when tracking the forward evolution of a random process, similarly $\bJ$ relates to a situation where we record new values of a random process at widely separated gaps of an earlier record. Regarding the significance of MR-spectra in time-series analysis, examples similar to the one that we presented here suggest similar qualities to the ME-ones (though, admitedly, they are slightly less appealing in terms of their ease of computation).

\end{document}